\documentclass[a4paper,12pt]{article}
    \usepackage[top=2.5cm,bottom=2.5cm,left=2.5cm,right=2.5cm]{geometry}
    \usepackage{cite, amsmath, amssymb}
    \usepackage[margin=1cm,%
                font=small,%
                format=hang,%
                labelsep=period,%
                labelfont=bf]{caption}
\usepackage[latin1]{inputenc}
\usepackage{rotating}
\usepackage{color}
\usepackage{graphicx}
\usepackage{epstopdf}
\usepackage{multirow} 
\usepackage{subfigure}
\usepackage{amsmath}
\usepackage{amsfonts}
\usepackage{latexsym}
\usepackage{authblk}
\usepackage[none]{hyphenat}
\begin{document}

\sloppy
\newcommand{\hide}[1]{}
\newcommand{\tbox}[1]{\mbox{\tiny #1}}
\newcommand{\half}{\mbox{\small $\frac{1}{2}$}}
\newcommand{\sinc}{\mbox{sinc}}
\newcommand{\const}{\mbox{const}}
\newcommand{\trc}{\mbox{trace}}
\newcommand{\intt}{\int\!\!\!\!\int }
\newcommand{\ointt}{\int\!\!\!\!\int\!\!\!\!\!\circ\ }
\newcommand{\eexp}{\mbox{e}^}
\newcommand{\bra}{\left\langle}
\newcommand{\ket}{\right\rangle}
\newcommand{\EPS} {\mbox{\LARGE $\epsilon$}}
\newcommand{\ar}{\mathsf r}
\newcommand{\im}{\mbox{Im}}
\newcommand{\re}{\mbox{Re}}
\newcommand{\bmsf}[1]{\bm{\mathsf{#1}}}
\newcommand{\mpg}[2][1.0\hsize]{\begin{minipage}[b]{#1}{#2}\end{minipage}}
\newcommand{\CC}{\mathbb{C}}
\newcommand{\NN}{\mathbb{N}}
\newcommand{\PP}{\mathbb{P}}
\newcommand{\RR}{\mathbb{R}}
\newcommand{\QQ}{\mathbb{Q}}
\newcommand{\ZZ}{\mathbb{Z}}
\renewcommand{\a}{\alpha}
\renewcommand{\b}{\beta}
\renewcommand{\d}{\delta}
\newcommand{\D}{\Delta}
\newcommand{\g}{\gamma}
\newcommand{\G}{\Gamma}
\renewcommand{\th}{\theta}
\renewcommand{\l}{\lambda}
\renewcommand{\L}{\Lambda}
\renewcommand{\O}{\Omega}
\newcommand{\s}{\sigma}

\newtheorem{theorem}{Theorem}
\newtheorem{acknowledgement}[theorem]{Acknowledgement}
\newtheorem{algorithm}[theorem]{Algorithm}
\newtheorem{axiom}[theorem]{Axiom}
\newtheorem{claim}[theorem]{Claim}
\newtheorem{conclusion}[theorem]{Conclusion}
\newtheorem{condition}[theorem]{Condition}
\newtheorem{conjecture}[theorem]{Conjecture}
\newtheorem{corollary}[theorem]{Corollary}
\newtheorem{criterion}[theorem]{Criterion}
\newtheorem{definition}[theorem]{Definition}
\newtheorem{example}[theorem]{Example}
\newtheorem{exercise}[theorem]{Exercise}
\newtheorem{lemma}[theorem]{Lemma}
\newtheorem{notation}[theorem]{Notation}
\newtheorem{problem}[theorem]{Problem}
\newtheorem{proposition}[theorem]{Proposition}
\newtheorem{remark}[theorem]{Remark}
\newtheorem{solution}[theorem]{Solution}
\newtheorem{summary}[theorem]{Summary}
\newenvironment{proof}[1][Proof]{\noindent\textbf{#1.} }{\ \rule{0.5em}{0.5em}}

\title{\vspace*{3.5cm}\textbf{Normalized Sombor indices as complexity measures of random graphs}}

\author[1]{\bf{R. Aguilar-S\'anchez}}
\author[2]{\bf{J. A. M\'endez-Berm\'udez}}
\author[3]{\bf{Jos\'e M. Rodr\'{\i}guez}}
\author[4]{\bf{Jos\'e M. Sigarreta\footnote{Corresponding author}}}

\affil[1]{\small{\it Facultad de Ciencias Qu\'imicas, Benem\'erita Universidad
Aut\'onoma de Puebla, Puebla 72570, Mexico}}
\affil[2]{\it Instituto de F\'{\i}sica, Benem\'erita Universidad Aut\'onoma de Puebla,
Apartado Postal J-48, Puebla 72570, Mexico}
\affil[3]{\it Departamento de Matem\'aticas, Universidad Carlos III de Madrid, 
Avenida de la Universidad 30, 28911 Legan\'es, Madrid, Spain}
\affil[4]{\it Facultad de Matem\'aticas, Universidad Aut\'onoma de Guerrero, 
Carlos E. Adame No.54 Col. Garita, Acapulco Gro. 39650, Mexico}
\affil[ ]{\ttfamily {\textbf {ras747698@gmail.com, jmendezb@ifuap.buap.mx, jomaro@math.uc3m.es, jsmathguerrero@gmail.com}}}
\date{}

\maketitle
\thispagestyle{empty}

\centerline{(Received xxx)}
\begin{abstract}
We perform a detailed computational study of the recently introduced Sombor indices on random 
graphs. Specifically, we apply Sombor indices on three models of random graphs: Erd\"os-R\'enyi graphs, 
random geometric graphs, and bipartite random graphs.
Within a statistical random matrix theory approach, we show that the average values 
of Sombor indices, normalized to the order of the graph, scale with the graph average degree. 
Moreover, we discuss the application of average Sombor indices as complexity measures of random 
graphs and, as a consequence, we show that selected normalized Sombor indices are highly 
correlated with the Shannon entropy of the eigenvectors of the graph adjacency matrix.
\end{abstract}

\baselineskip=0.30in

\section{Introduction}

Given a graph $G=(V(G),E(G))$, the Sombor index of $G$, introduced by I. Gutman in \cite{G21a}, is defined as
\begin{equation}
SO(G)
= \sum_{uv\in E(G)} \sqrt{k_u^2 + k_v^2} ,
\label{SO}
\end{equation}
where $uv$ denotes the edge of the graph $G$ connecting the vertices $u$ and $v$ and $k_u$ is the 
degree of the vertex $u$. Also, the modified Sombor index of $G$ was proposed in \cite{KG21} as
\begin{equation}
^mSO(G)
= \sum_{uv\in E(G)} \frac{1}{\sqrt{k_u^2 + k_v^2}} .
\label{mSO}
\end{equation}
In addition, two other Sombor indices have been introduced:
the first Banhatti-Sombor index \cite{LZKM21}
\begin{equation}
BSO(G)
= \sum_{uv\in E(G)} \sqrt{\frac{1}{k_u^2} + \frac{1}{k_v^2} }
\label{BSO}
\end{equation}
and the $\alpha$-Sombor index \cite{RDA21}
\begin{equation}
SO_\alpha(G)
= \sum_{uv\in E(G)} (k_u^\alpha + k_v^\alpha)^{1/\alpha} ,
\label{pSO}
\end{equation}
here $\alpha \in \RR$.
In fact, there is a general index that includes all the Sombor indices listed above: the first $(\alpha,\beta)-KA$ index of $G$ 
which was introduced in \cite{K19} as
\begin{equation}
KA^1_{\alpha,\beta}(G)
= \sum_{uv\in E(G)} \left( k_u^\alpha + k_v^\alpha \right)^\beta ,
\label{KA1}
\end{equation}
with $\alpha,\beta \in \RR$.
Note that $SO(G)=KA^1_{2,1/2}(G)$, $^mSO(G)=KA^1_{2,-1/2}(G)$, $BSO(G)=KA^1_{-2,1/2}(G)$, and 
$SO_\alpha(G)=KA^1_{\alpha,1/\alpha}(G)$. Also, we note that $KA^1_{1,\beta}(G)$ equals the general 
sum-connectivity index \cite{ZT10}
$
\chi_\beta(G) = \sum_{uv \in E(G)} (k_u + k_v)^\beta .
$

Reduced versions of $SO(G)$, $^mSO(G)$ and $KA^1_{\alpha,\beta}(G)$ were also introduced 
in \cite{G21a,KG21,K21}.
However, when dealing with random graphs we use to approximate vertex degrees by average degrees 
and since average degrees may be less than one, reduced degree-based indices are not amenable for us.
Thus we do not consider reduced Sombor indices here.

Even though Sombor indices were introduced very recently, there are already several works available 
in the literature where these indices are applied to chemical graphs of interest, see 
e.g. \cite{RDA21,K19,K21,CGR21,CR21,GA21,AG21,FYL21,DTW21,L21a,L21b,LYH22,ZLM21a,ZLM21b}.
Also, bounds for Sombor indices as well as relations among them and with many other topological 
indices have been reported 
in \cite{RDA21,GA21,G21b,MMM21,DCCS21,RRS21,MMAM21,WMLF21}.
From the application point of view, they have been shown to be useful to to model entropy and enthalpy of 
vaporization of alkanes \cite{R21}. 
In addition, the Sombor matrix has been proposed and studied in \cite{LM21}.
However, to the best of our knowledge, Sombor indices have not been applied to random graphs yet; 
thus in this work we undertake this task.

\bigskip

Here we consider three models of random graphs: Erd\"os-R\'enyi (ER) graphs, random geometric (RG) 
graphs, and bipartite random (BR) graphs. 
ER graphs~\cite{SR51,ER59,ER60} $G_{\tbox{ER}}(n,p)$ are formed by $n$ vertices connected independently 
with probability $p \in [0,1]$. 
While RG graphs~\cite{DC02,P03} $G_{\tbox{RG}}(n,r)$ consist of $n$ vertices uniformly and independently 
distributed on the unit square, where two vertices are connected by an edge if their Euclidean distance is less 
or equal than the connection radius $r \in [0,\sqrt{2}]$.
In addition we examine BR graphs $G_{\tbox{BR}}(n_1,n_2,p)$ composed by two disjoint sets, set 1 and set 2, 
with $n_1$ and $n_2$ vertices each such that there are no adjacent vertices within the same set, being 
$n=n_1+n_2$ the total number of vertices in the bipartite graph. The vertices of the two sets are connected 
randomly with probability $p \in [0,1]$.

We stress that the computational study of Sombor indices we perform here is justified by the random nature 
of the graph models we want to explore. Since a given parameter set [$(n,p)$, $(n,r)$, or $(n_1,n_2,p)$] 
represents an infinite-size ensemble of random [ER, RG, or BR] graphs, the computation of a Sombor index 
on a single graph is irrelevant. In contrast, the computation of a Sombor index on a large ensemble of random 
graphs, all characterized by the same parameter set, may provide useful {\it average} information about the 
full ensemble. This {\it statistical} approach, well known in random matrix theory studies, is not widespread in 
studies involving topological indices, mainly because topological indices are not commonly applied to random 
networks; for very recent exceptions see~\cite{MMRS20,MMRS21}.

Therefore, the purpose of this work is threefold. First, we push forward the statistical (computational) analysis 
of topological indices as a generic tool for studying average properties of random graphs; second, we perform 
for the first time (to our knowledge), a scaling study of Sombor indices on random graphs; and third, we discuss 
the application of selected Sombor indices as complexity measures of random graphs.

\section{Computational properties of Sombor indices on random graphs}

\subsection{Sombor indices on Erd\"os-R\'enyi graphs}
\label{ER}

In what follows we present the average values of the indices defined in Eqs.~(\ref{SO}-\ref{KA1}). 
All averages are computed over ensembles of $10^7/n$ ER graphs characterized 
by the parameter pair $(n,p)$.

\begin{figure}[t!]
\begin{center} 
\includegraphics[width=0.75\textwidth]{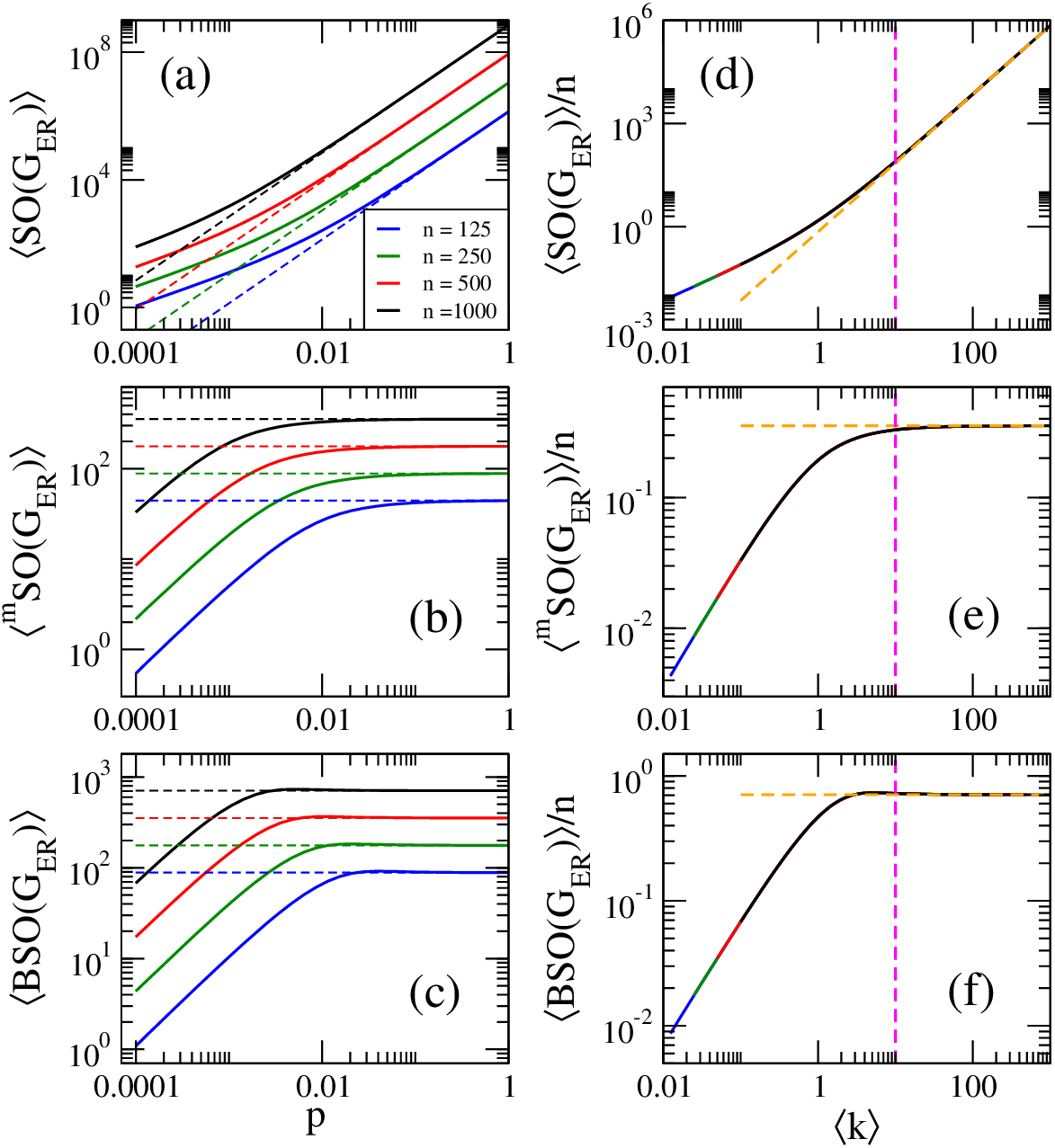}
\caption{\footnotesize{
(a) Average Sombor index $\left< SO(G_{\tbox{ER}}) \right>$,
(b) average modified Sombor index $\left< ^mSO(G_{\tbox{ER}}) \right>$, and
(c) average first Banhatti-Sombor index $\left< BSO(G_{\tbox{ER}}) \right>$ as a 
function of the probability $p$ of Erd\"os-R\'enyi graphs of size $n$.
(d) $\left< SO(G_{\tbox{ER}}) \right>/n$,
(e) $\left< ^mSO(G_{\tbox{ER}}) \right>/n$, and 
(f) $\left< BSO(G_{\tbox{ER}}) \right>/n$ as a function of $\left< k \right>$.
Dashed lines in panels (a), (b) and (c) correspond to Eqs.~(\ref{SOp}), (\ref{mSOp}) and (\ref{BSOp}), respectively.
While dashed lines in panels (d), (e) and (f) are Eqs.~(\ref{SOk}), (\ref{mSOk}) and (\ref{BSOk}), respectively.
The vertical magenta dashed line in (b-f) marks $\left< k \right>= 10$.
}}
\label{Fig01}
\end{center}
\end{figure}

\begin{figure}[t!]
\begin{center} 
\includegraphics[width=0.75\textwidth]{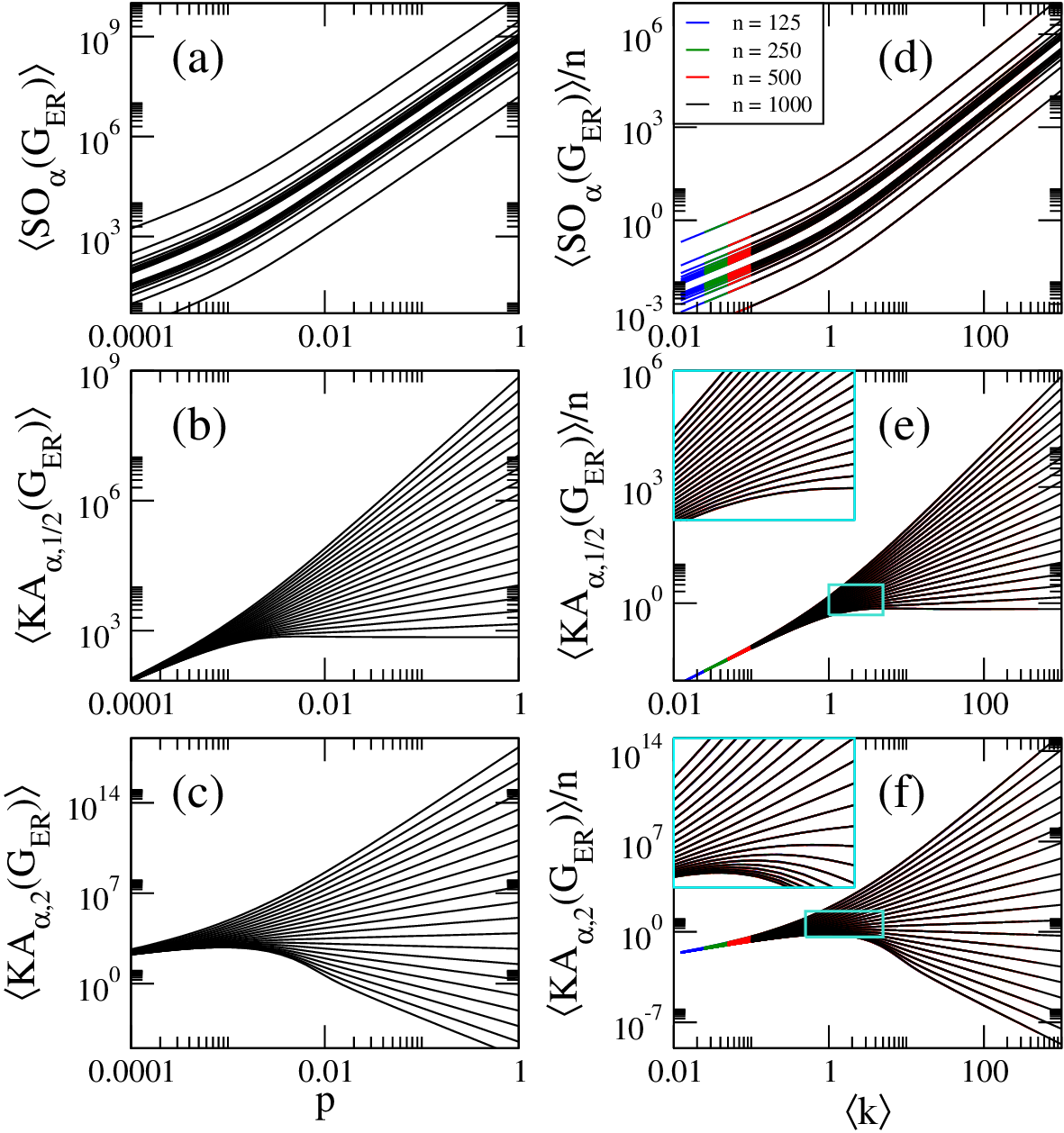}
\caption{\footnotesize{
(a) Average $\alpha$-Sombor index $\left< SO_\alpha(G_{\tbox{ER}}) \right>$,
(b) average first $(\alpha,\beta)-KA$ index $\left< KA_{\alpha,\beta}(G_{\tbox{ER}}) \right>$, with $\beta=1/2$, and
(c) average first $(\alpha,\beta)-KA$ index $\left< KA_{\alpha,\beta}(G_{\tbox{ER}}) \right>$, with $\beta=2$, 
as a function of the probability $p$ of Erd\"os-R\'enyi graphs of size $n=1000$.
In all panels we show curves for $\alpha\in[-2,2]$ in steps of 0.2 (from bottom to top).
(d) $\left< SO_\alpha(G_{\tbox{ER}}) \right>/n$,
(e) $\left< KA_{\alpha,1/2}(G_{\tbox{ER}}) \right>/n$, and
(f) $\left< KA_{\alpha,2}(G_{\tbox{ER}}) \right>/n$ as a function of $\left< k \right>$
for ER graphs of four different sizes $n$.
The insets in (e,f) are enlargements of the cyan rectangles.
}}
\label{Fig02}
\end{center}
\end{figure}

On the one hand, in Figs.~\ref{Fig01}(a), \ref{Fig01}(b), and~\ref{Fig01}(c) we present, respectively, 
the average Sombor index $\left< SO(G_{\tbox{ER}}) \right>$, 
the average modified Sombor index $\left< ^mSO(G_{\tbox{ER}}) \right>$, and
the average first Banhatti-Sombor index $\left< BSO(G_{\tbox{ER}}) \right>$
as a function of the probability $p$ of ER graphs of sizes $n=\{125,250,500,1000\}$.
On the other hand, in Fig.~\ref{Fig02} we plot
the average $\alpha$-Sombor index $\left< SO_\alpha(G_{\tbox{ER}}) \right>$, see Fig.~\ref{Fig02}(a),
and the average first $(\alpha,\beta)-KA$ index $\left< KA^1_{\alpha,\beta}(G_{\tbox{ER}}) \right>$, see 
Figs.~\ref{Fig02}(c,d), as a function $p$ of ER graphs of size $n=1000$.
In Fig.~\ref{Fig02} we show curves for $\alpha\in[-2,2]$ and, in the case of 
$\left< KA^1_{\alpha,\beta}(G_{\tbox{ER}}) \right>$,
we choose to report $\beta=1/2$ and $\beta=2$ as representative cases.

From this figures we observe that:
\begin{itemize}

\item[{\bf (i)}]
The curves of $\left< SO(G_{\tbox{ER}}) \right>$ and $\left< SO_\alpha(G_{\tbox{ER}}) \right>$
are monotonically increasing functions of $p$. See Figs.~\ref{Fig01}(a) and~\ref{Fig02}(a).

\item[{\bf (ii)}]
The curves of $\left< ^mSO(G_{\tbox{ER}}) \right>$ and $\left< BSO(G_{\tbox{ER}}) \right>$
grow for small $p$ and saturate above a given value of $p$.
See Figs.~\ref{Fig01}(b) and~\ref{Fig01}(c).

\item[{\bf (iii)}]
The curves of $\left< KA^1_{\alpha,\beta}(G_{\tbox{ER}}) \right>$ show three
different behaviors as a function of $p$ depending on the values of $\alpha$ and $\beta$:
For $\alpha<\alpha_0$, they grow for small $p$, approach a 
maximum value and then decrease when $p$ is further increased. 
For $\alpha>\alpha_0$, they are monotonically increasing functions of $p$.
For $\alpha=\alpha_0$ the curves saturate above a given value of $p$.
For $\beta=1/2$ and $\beta=2$, the cases reported in Figs.~\ref{Fig02}(c,d), we found
$\alpha_0=-2$ and $\alpha_0=-1/2$, respectively.

\item[{\bf (iv)}]
When $np\gg 1$, we can approximate $k_u \approx k_v \approx \left<  k \right>$
in Eqs.~(\ref{SO}-\ref{KA1}), with
\begin{equation}
\label{k}
\left< k \right> = (n-1)p .
\end{equation}
Therefore, for $np\gg 1$, the average values of the Sombor indices are well approximated by:
\begin{equation}
\label{SOp}
\left< SO(G_{\tbox{ER}}) \right> \approx  \frac{n}{\sqrt{2}} \left[ (n-1)p \right]^2 ,
\end{equation}
\begin{equation}
\label{mSOp}
\left< ^mSO(G_{\tbox{ER}}) \right> \approx \frac{n}{2\sqrt{2}} ,
\end{equation}
\begin{equation}
\label{BSOp}
\left< BSO(G_{\tbox{ER}}) \right> \approx  \frac{n}{\sqrt{2}} ,
\end{equation}
\begin{equation}
\label{SOap}
\left< SO_\alpha(G_{\tbox{ER}}) \right> \approx \frac{n}{2^{1-1/\alpha}} \left[ (n-1)p \right]^2 ,
\end{equation}
\begin{equation}
\label{KAabp}
\left< KA^1_{\alpha,\beta}(G_{\tbox{ER}}) \right> \approx \frac{n}{2^{1-\beta}} \left[ (n-1)p \right]^{1+\alpha\beta} .
\end{equation}
In Figs.~\ref{Fig01}(a)-\ref{Fig01}(c), we show that Eqs.~(\ref{SOp}-\ref{BSOp}) (dashed 
lines) indeed describe well the data (thick full curves) for large enough $p$.
We also verified that Eqs.~(\ref{SOap},\ref{KAabp}) describe well the data for $np\gg 1$ reported in Figs.~\ref{Fig02}(a-c),
however we did not include them to avoid figure saturation.

\end{itemize}


We note that in Figs.~\ref{Fig01}(a-c) we present average Sombor indices as a function of the 
probability $p$ of ER graphs of four different sizes $n$. 
It is quite clear from these figures that the curves, characterized by the different network sizes, 
are very similar but displaced on both axes. 
A similar observation can be made for $\left< SO_\alpha(G_{\tbox{ER}}) \right>$ and 
$\left< KA_{\alpha,\beta}(G_{\tbox{ER}}) \right>$ (not shown in Figs.~\ref{Fig02}(a-c) to avoid 
figure saturation).
This behavior suggests that the average Sombor indices can be scaled. Then, in what follows
we look for the parameters that scale the average Sombor indices.

From Eqs.~(\ref{SOp}-\ref{KAabp}) we observe that $\left< X(G_{\tbox{ER}}) \right>\propto n f[(n-1)p)]$ or
\begin{equation}
\label{scaling}
\left< X(G_{\tbox{ER}}) \right>\propto n f(\left< k \right>),
\end{equation}
where $X$ represents all the Sombor indices studied here.
Therefore, in Figs.~\ref{Fig01}(d-f) and~\ref{Fig02}(d-f) we plot average Sombor indices,
normalized to $n$, as a function of $\left< k \right>$ showing that all indices are now properly scaled; 
i.e.~the curves painted in different colors for different graph sizes fall on top of each other. 
Moreover, we can rewrite Eqs.~(\ref{SOp}-\ref{KAabp}) as
\begin{equation}
\label{SOk}
\frac{\left< SO(G_{\tbox{ER}}) \right>}{n} \approx \frac{1}{\sqrt{2}} \left< k \right>^2 ,
\end{equation}
\begin{equation}
\label{mSOk}
\frac{\left< ^mSO(G_{\tbox{ER}}) \right>}{n} \approx \frac{1}{2\sqrt{2}} ,
\end{equation}
\begin{equation}
\label{BSOk}
\frac{\left< BSO(G_{\tbox{ER}}) \right>}{n} \approx \frac{1}{\sqrt{2}} ,
\end{equation}
\begin{equation}
\label{SOak}
\frac{\left< SO_\alpha(G_{\tbox{ER}}) \right>}{n} \approx \frac{1}{2^{1-1/\alpha}} \left< k \right>^2 ,
\end{equation}
\begin{equation}
\label{KAabk}
\frac{\left< KA^1_{\alpha,\beta}(G_{\tbox{ER}}) \right>}{n} \approx \frac{1}{2^{1-\beta}} \left< k \right>^{1+\alpha\beta} .
\end{equation}
In Figs.~\ref{Fig01}(d)-\ref{Fig01}(f), we show that Eqs.~(\ref{SOk}-\ref{BSOk}) (orange-dashed 
lines) indeed describe well the data (thick full curves) for $\left<  k \right>\ge 10$.
We also verified that Eqs.~(\ref{SOak}-\ref{KAabk}) describe well the data for $\left<  k \right>\ge 10$ reported in 
Figs.~\ref{Fig02}(d)-\ref{Fig02}(f) (not shown here to avoid figure saturation).

It is relevant to stress that even when Eq.~(\ref{scaling}) was deduced form Eqs.~(\ref{SOp}-\ref{KAabp}),
expected to be valid in the dense limit (i.e.~for $\left<  k \right> \gg 1$), it is indeed valid for any $\left<  k \right>$
as clearly seen in Figs.~\ref{Fig01}(d)-\ref{Fig01}(f) and Figs.~\ref{Fig02}(d)-\ref{Fig02}(f).

\subsection{Sombor indices on random geometric graphs}
\label{RG}

As in the previous Subsection, here we present the average values of the Sombor indices listed in 
Eqs.~(\ref{SO}-\ref{KA1}). Again, all averages are computed over ensembles of $10^7/n$ 
random graphs, each ensemble characterized by a fixed parameter pair $(n,r)$.

\begin{figure}[t!]
\begin{center} 
\includegraphics[width=0.75\textwidth]{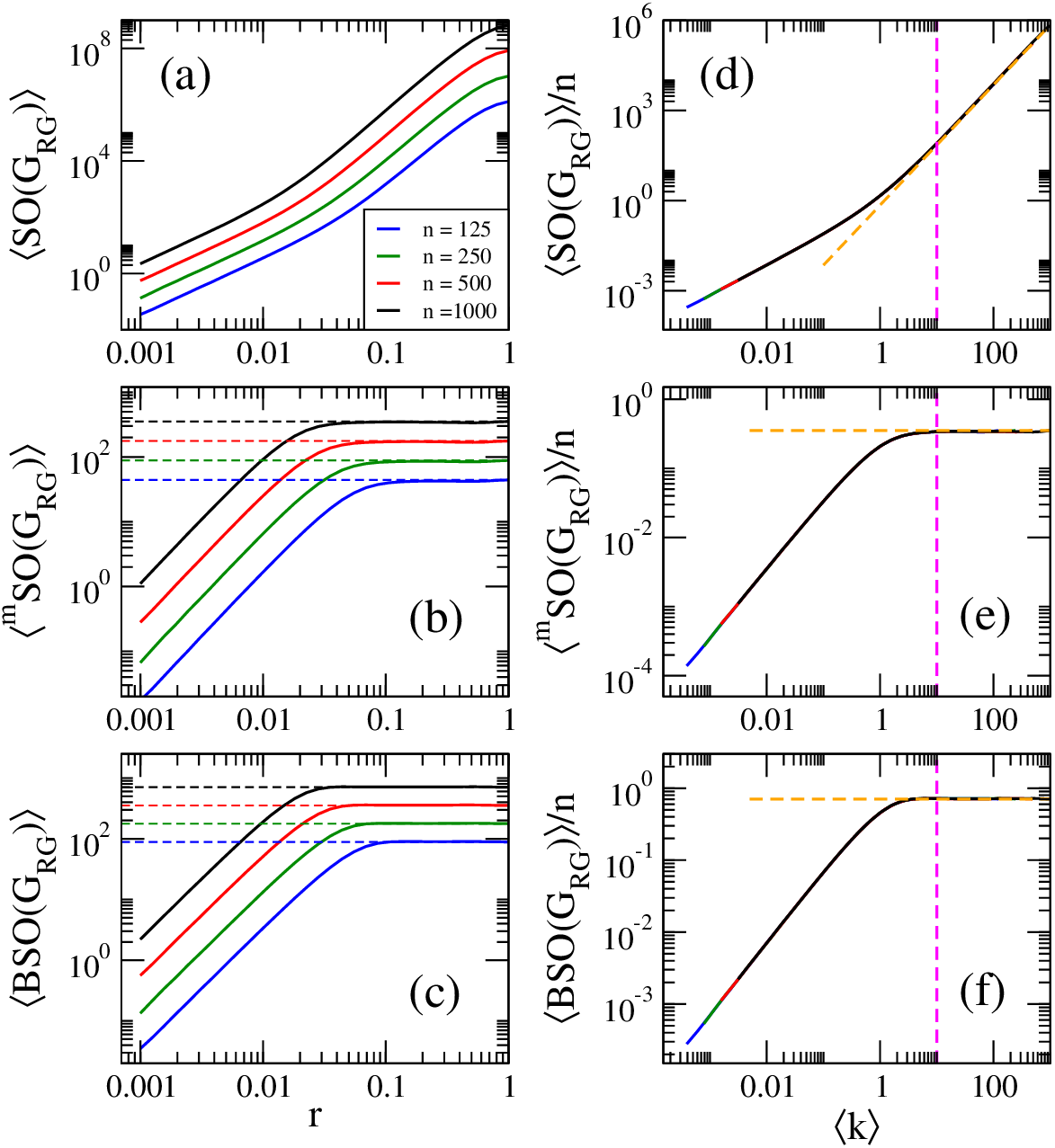}
\caption{\footnotesize{
(a) Average Sombor index $\left< SO(G_{\tbox{RG}}) \right>$,
(b) average modified Sombor index $\left< ^mSO(G_{\tbox{RG}}) \right>$, and
(c) average first Banhatti-Sombor index $\left< BSO(G_{\tbox{RG}}) \right>$ as a 
function of the connection radius $r$ of random geometric graphs of size $n$.
(d) $\left< SO(G_{\tbox{RG}}) \right>/n$,
(e) $\left< ^mSO(G_{\tbox{RG}}) \right>/n$, and 
(f) $\left< BSO(G_{\tbox{RG}}) \right>/n$ as a function of $\left< k \right>$.
Dashed lines in panels (a), (b) and (c) correspond to Eqs.~(\ref{SOr}), (\ref{mSOr}) and (\ref{BSOr}), respectively.
While dashed lines in panels (d), (e) and (f) are Eqs.~(\ref{SOk}), (\ref{mSOk}) and (\ref{BSOk}), respectively.
The vertical magenta dashed line in (b-f) marks $\left< k \right>= 10$.
}}
\label{Fig03}
\end{center}
\end{figure}
\begin{figure}[t!]
\begin{center} 
\includegraphics[width=0.75\textwidth]{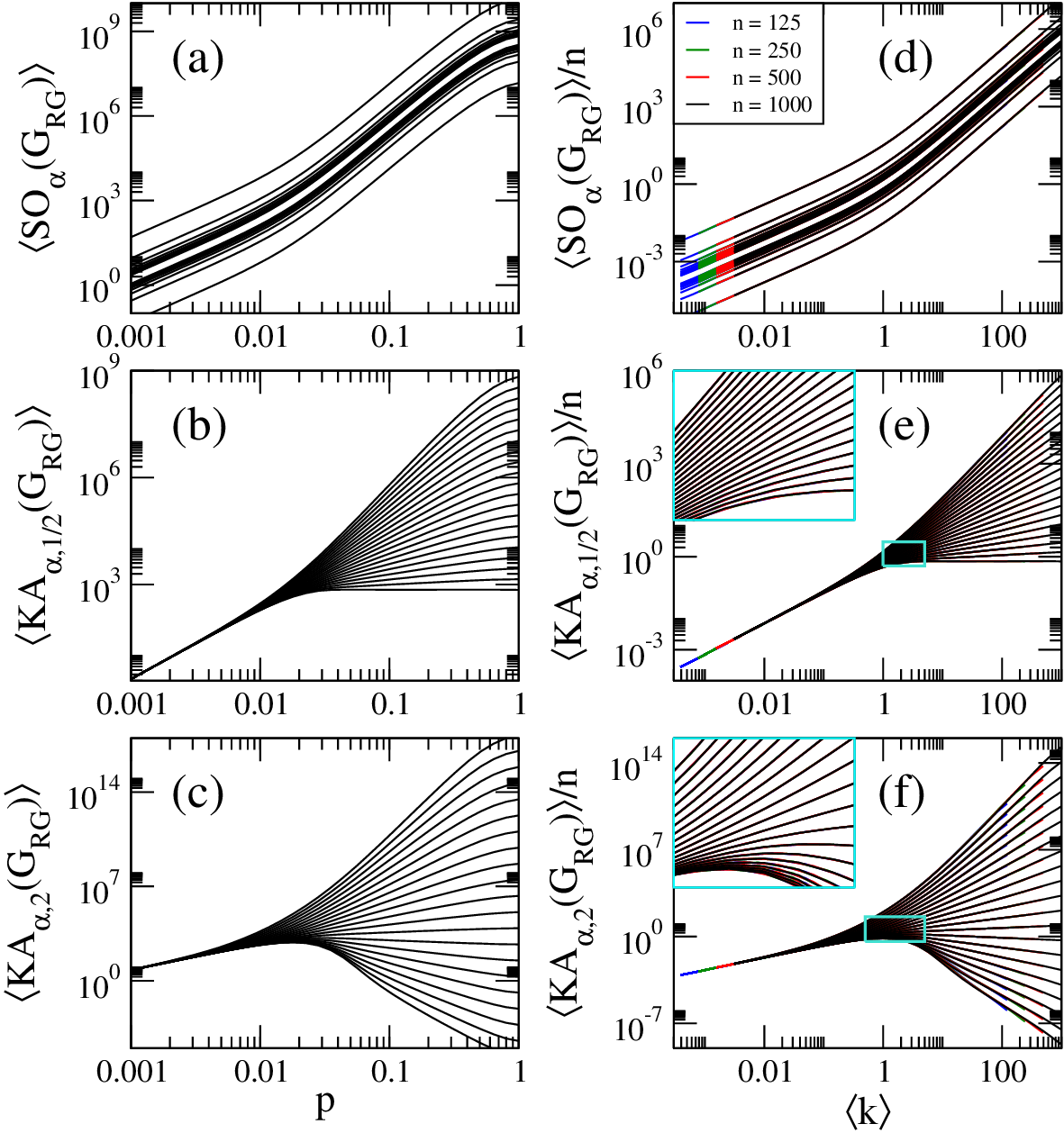}
\caption{\footnotesize{
(a) Average $\alpha$-Sombor index $\left< SO_\alpha(G_{\tbox{RG}}) \right>$,
(b) average first $(\alpha,\beta)-KA$ index $\left< KA_{\alpha,\beta}(G_{\tbox{RG}}) \right>$, with $\beta=1/2$, and
(c) average first $(\alpha,\beta)-KA$ index $\left< KA_{\alpha,\beta}(G_{\tbox{RG}}) \right>$, with $\beta=2$, 
as a function of the connection radius $r$ of random geometric graphs of size $n=1000$.
In all panels we show curves for $\alpha\in[-2,2]$ in steps of 0.2 (from bottom to top).
(d) $\left< SO_\alpha(G_{\tbox{RG}}) \right>/n$,
(e) $\left< KA_{\alpha,1/2}(G_{\tbox{RG}}) \right>/n$, and
(f) $\left< KA_{\alpha,2}(G_{\tbox{RG}}) \right>/n$ as a function of $\left< k \right>$
for RG graphs of four different sizes $n$.
The insets in (e,f) are enlargements of the cyan rectangles.
}}
\label{Fig04}
\end{center}
\end{figure}

Then, in Figs.~\ref{Fig03}(a), \ref{Fig03}(b), and~\ref{Fig03}(c) we present, respectively, 
the average Sombor index $\left< SO(G_{\tbox{RG}}) \right>$, 
the average modified Sombor index $\left< ^mSO(G_{\tbox{RG}}) \right>$, and
the average first Banhatti-Sombor index $\left< BSO(G_{\tbox{RG}}) \right>$
as a function of the connection radius $r$ of RG graphs of sizes $n=\{125,250,500,1000\}$.
Also, in Fig.~\ref{Fig04} we plot
the average $\alpha$-Sombor index $\left< SO_\alpha(G_{\tbox{RG}}) \right>$, see Fig.~\ref{Fig04}(a),
and the average first $(\alpha,\beta)-KA$ index $\left< KA^1_{\alpha,\beta}(G_{\tbox{RG}}) \right>$, see 
Figs.~\ref{Fig04}(c,d), as a function $r$ of RG graphs of size $n=1000$.

For comparison purposes, Figs.~\ref{Fig03} and~\ref{Fig04} are similar to 
Figs.~\ref{Fig01} and~\ref{Fig02}. In fact, all the observations {\bf (i-iv)} made in the previous 
Subsection for ER graphs are also valid for RG graphs by replacing $G_{\tbox{ER}}\to G_{\tbox{RG}}$ 
and $p\to g(r)$, with~\cite{EM15}
\begin{equation}
g(r) = 
\left\{ 
\begin{array}{ll}
           r^2  \left[ \pi - \frac{8}{3}r +\frac{1}{2}r^2 \right] & 0 \leq r \leq 1 \, , 
           \vspace{0.25cm} \\
             \frac{1}{3} - 2r^2 \left[ 1 - \arcsin(1/r) + \arccos(1/r) \right] 
             +\frac{4}{3}(2r^2+1) \sqrt{r^2-1} 
             -\frac{1}{2}r^4 & 1 \leq r \leq \sqrt{2} \, .
\end{array}
\right.
\label{g(r)}
\end{equation}
However, given the fact that this is the first study (to our knowledge) of average Sombor indices on RG 
graphs, we want to stress that when $nr\gg 1$, we can approximate $k_u \approx k_v \approx \left<  k \right>$
in Eqs.~(\ref{SO}-\ref{KA1}), with
\begin{equation}
\label{kRG}
\left< k \right> = (n-1)g(r) .
\end{equation}
Therefore, in the dense limit, the average values of the Sombor indices on RG graphs are well approximated by:
\begin{equation}
\label{SOr}
\left< SO(G_{\tbox{RG}}) \right> \approx  \frac{n}{\sqrt{2}} \left[ (n-1)g(r) \right]^2 ,
\end{equation}
\begin{equation}
\label{mSOr}
\left< ^mSO(G_{\tbox{RG}}) \right> \approx \frac{n}{2\sqrt{2}} ,
\end{equation}
\begin{equation}
\label{BSOr}
\left< BSO(G_{\tbox{RG}}) \right> \approx  \frac{n}{\sqrt{2}} ,
\end{equation}
\begin{equation}
\label{SOar}
\left< SO_\alpha(G_{\tbox{RG}}) \right> \approx \frac{n}{2^{1-1/\alpha}} \left[ (n-1)g(r) \right]^2 ,
\end{equation}
\begin{equation}
\label{KAabr}
\left< KA^1_{\alpha,\beta}(G_{\tbox{RG}}) \right> \approx \frac{n}{2^{1-\beta}} \left[ (n-1)g(r) \right]^{1+\alpha\beta} .
\end{equation}
In Figs.~\ref{Fig03}(a)-\ref{Fig03}(c), we show that Eqs.~(\ref{SOr}-\ref{BSOr}) (dashed 
lines) indeed describe well the data (thick full curves) for large enough $r$.
We also verified that Eqs.~(\ref{SOar},\ref{KAabr}) describe well the data reported in Figs.~\ref{Fig04}(a-c),
for large enough $r$,
however we did not include them to avoid figure saturation.

It is quite remarkable to note that by substituting the average degree of Eq.~(\ref{kRG}) into 
Eqs.~(\ref{SOr}-\ref{BSOr}) we get exactly the same expressions listed in Eqs.~(\ref{SOk}-\ref{KAabk}).
Therefore, in Figs.~\ref{Fig03}(d-f) and~\ref{Fig04}(d-f) we plot average Sombor indices,
on RG graphs, normalized to $n$, as a function of $\left< k \right>$ showing that all curves are properly scaled.
Also, in Figs.~\ref{Fig03}(d)-\ref{Fig03}(f), we show that Eqs.~(\ref{SOk}-\ref{BSOk}) (orange-dashed 
lines) indeed describe well the data (thick full curves) for $\left<  k \right>\ge 10$.
We also verified (not shown here) that Eqs.~(\ref{SOak}-\ref{KAabk}) describe well the data for $\left<  k \right>\ge 10$ 
reported in Figs.~\ref{Fig02}(d)-\ref{Fig02}(f).

\subsection{Sombor indices on bipartite random graphs}
\label{BR}

Now we compute average Sombor indices on ensembles of $10^7/n$ BR graphs. In contrast to 
ER and RG graphs now the BR graph ensembles are characterized by three parameters: $n_1$, $n_2$, and $p$. 
Thus we consider two cases: $n_1=n_2$ and $n_1<n_2$.
We note that bounds for the Sombor index on bipartite graphs have been reported in \cite{DCCS21}.

\begin{figure}[t!]
\begin{center} 
\includegraphics[width=0.75\textwidth]{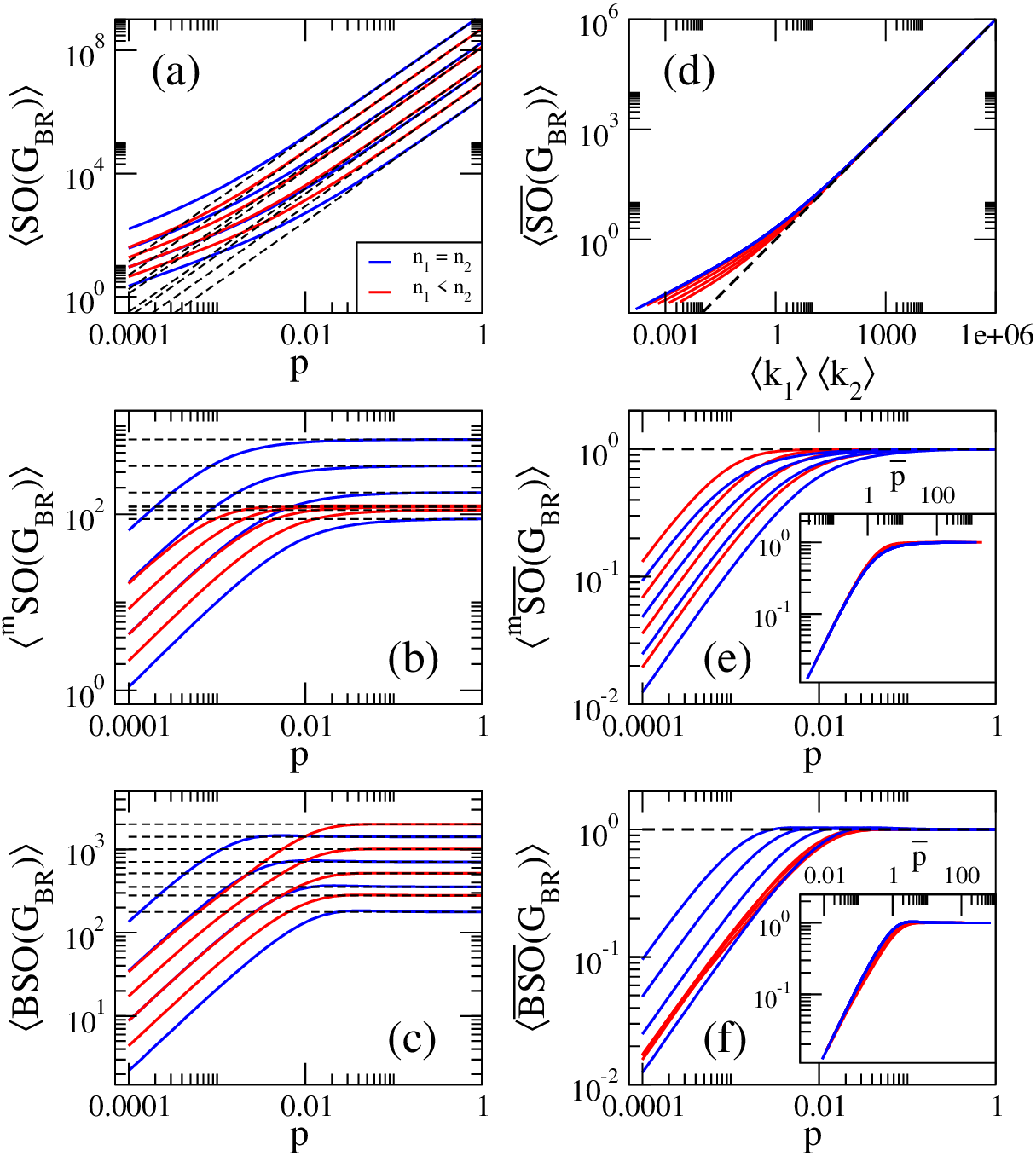}
\caption{\footnotesize{
(a) Average Sombor index $\left< SO(G_{\tbox{BR}}) \right>$,
(b) average modified Sombor index $\left< ^mSO(G_{\tbox{BR}}) \right>$, and
(c) average first Banhatti-Sombor index $\left< BSO(G_{\tbox{BR}}) \right>$ as a 
function of the probability $p$ of bipartite random graphs with sets of seizes $n_1$ and $n_2$.
In all panels: $n_1=n_2$ with $n_2=\{125,250,500,1000\}$ (blue lines, $n_2$ increases from bottom to top)
and $n_1<n_2$ with $n_1=125$ and $n_2=\{250,500,1000,2000\}$ (red lines, $n_2$ increases from bottom to top).
(d) $\left< \overline{SO}(G_{\tbox{BR}}) \right>=\left< SO(G_{\tbox{BR}}) \right>/\sqrt{n_1^2+n_2^2}$ 
vs.~the product $\left< k_1 \right>\left< k_2 \right>$.
(e) $\left< ^m\overline{SO}(G_{\tbox{BR}}) \right>=\sqrt{n_1^2+n_2^2} \left< ^mSO(G_{\tbox{BR}}) \right>/(n_1n_2)$ 
vs.~$p$.
(f) $\left< \overline{BSO}(G_{\tbox{BR}}) \right>=\left< BSO(G_{\tbox{BR}}) \right>/\sqrt{n_1^2+n_2^2}$ vs.~$p$.
Dashed lines in panels (a), (b) and (c) correspond to Eqs.~(\ref{SOBGp}), (\ref{mSOBGp}) and (\ref{BSOBGp}), respectively.
While dashed lines in panels (d), (e) and (f) are Eqs.~(\ref{SOBGk}), (\ref{mSOBGp}) and (\ref{BSOBGp}), respectively.
The inset in (e) shows $\left< ^m\overline{SO}(G_{\tbox{BR}}) \right>$ vs.~$\overline{p}=p\sqrt{n_1^2+n_2^2}$.
The inset in (f) shows $\left< \overline{BSO}(G_{\tbox{BR}}) \right>$ vs.~$\overline{p}=pn_1n_2/\sqrt{n_1^2+n_2^2}$.
}}
\label{Fig05}
\end{center}
\end{figure}

\begin{figure}[t!]
\begin{center} 
\includegraphics[width=0.75\textwidth]{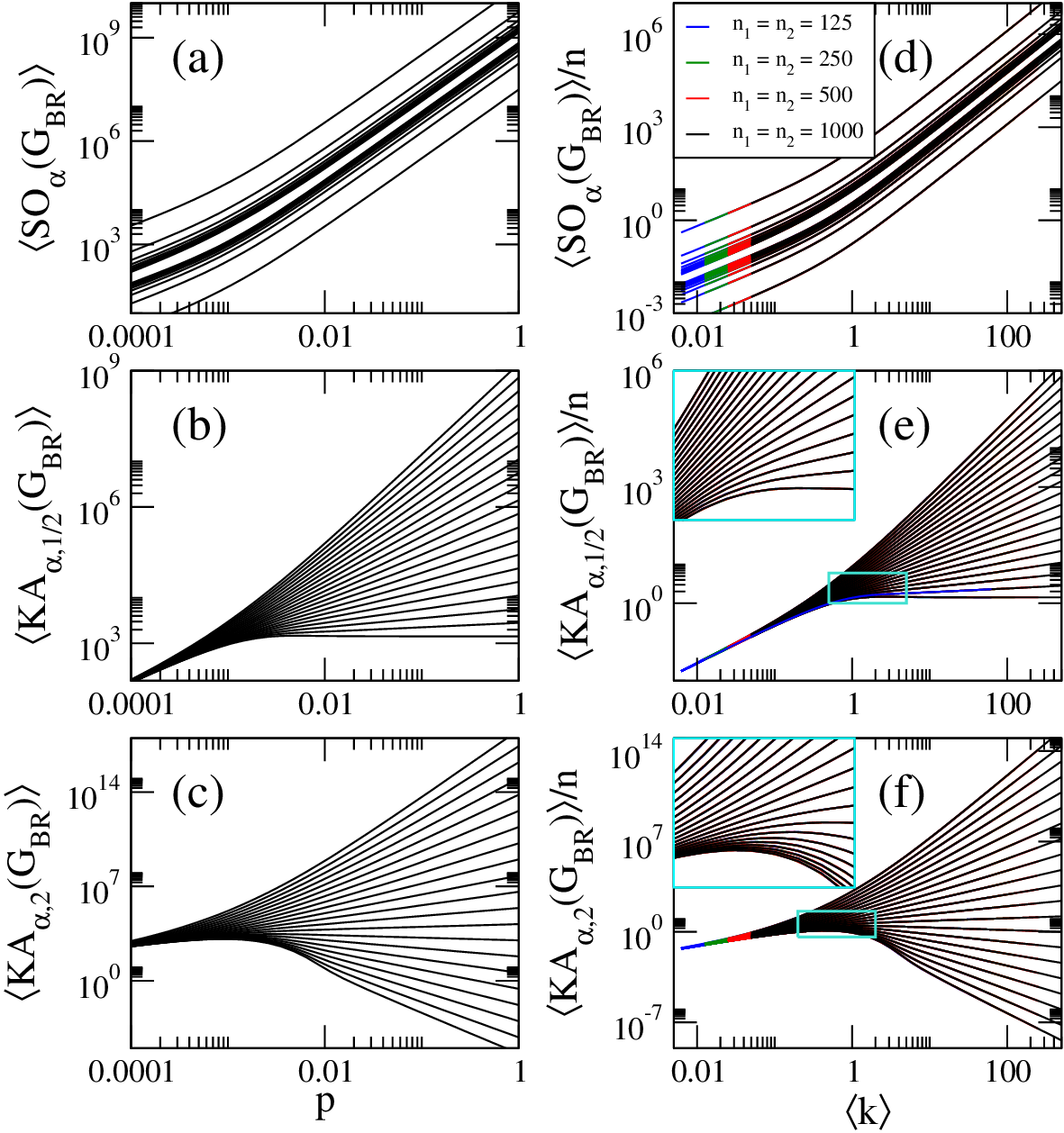}
\caption{\footnotesize{
(a) Average $\alpha$-Sombor index $\left< SO_\alpha(G_{\tbox{BR}}) \right>$, 
(b) average first $(\alpha,\beta)-KA$ index $\left< KA_{\alpha,\beta}(G_{\tbox{BR}}) \right>$, with $\beta=1/2$, and
(c) average first $(\alpha,\beta)-KA$ index $\left< KA_{\alpha,\beta}(G_{\tbox{BR}}) \right>$, with $\beta=2$, 
as a function of the probability $p$ of bipartite random graphs with sets of seizes $n_1=n_2=1000$.
In all panels we show curves for $\alpha\in[-2,2]$ in steps of 0.2 (from bottom to top).
(d) $\left< SO_\alpha(G_{\tbox{BR}}) \right>/n$,
(e) $\left< KA_{\alpha,1/2}(G_{\tbox{BR}}) \right>/n$, and
(f) $\left< KA_{\alpha,2}(G_{\tbox{BR}}) \right>/n$ as a function of $\left< k \right>$
for BR graphs of four different sizes $n$.
The insets in (e,f) are enlargements of the cyan rectangles.
}}
\label{Fig06}
\end{center}
\end{figure}

In Figs.~\ref{Fig05}(a), \ref{Fig05}(b), and~\ref{Fig05}(c) we present, respectively, 
the average Sombor index $\left< SO(G_{\tbox{BR}}) \right>$, 
the average modified Sombor index $\left< ^mSO(G_{\tbox{BR}}) \right>$, and
the average first Banhatti-Sombor index $\left< BSO(G_{\tbox{BR}}) \right>$
as a function of the probability $p$ of BR graphs characterized by $n_1=n_2$ with $n_2=\{125,250,500,1000\}$
(blue lines) and BR graphs characterized by $n_1<n_2$ with $n_1=125$ and $n_2=\{125,250,500,1000\}$ (red lines).
Also, in Fig.~\ref{Fig06} we plot
the average $\alpha$-Sombor index $\left< SO_\alpha(G_{\tbox{BR}}) \right>$, see Fig.~\ref{Fig06}(a),
and the average first $(\alpha,\beta)-KA$ index $\left< KA^1_{\alpha,\beta}(G_{\tbox{BR}}) \right>$, see 
Figs.~\ref{Fig06}(c,d), as a function $p$ of BR graphs of size $n_1=n_2=1000$.

It is interesting to notice that all the observations {\bf (i-iv)} made in Subsection~\ref{RG} for ER graphs are 
also valid for BR graphs by just replacing $G_{\tbox{ER}}\to G_{\tbox{BR}}$.
Moreover, we can also write approximate expressions for the average Sombor indices on BR graphs
in the dense limit.
However, since edges in a bipartite graph join vertices of different sets, and we are labeling 
here the sets as set 1 and set 2, we replace $d_u$ by $d_1$ and $d_v$ by $d_2$ in the expression 
for the Sombor indices. 
Thus, when $n_1p\gg 1$ and $n_2p\gg 1$, we can approximate $k_u = k_1 \approx \left<  k_1 \right>$
and $k_v = k_2 \approx \left<  k_2 \right>$ in Eqs.~(\ref{SO}-\ref{KA1}), with
\begin{equation}
\label{kBR}
\left< k_{1,2} \right> = n_{2,1}p .
\end{equation}
Therefore, in the dense limit, the average values of the Sombor indices on BR graphs are 
well approximated by:
\begin{equation}
\label{SOBGp}
\left< SO(G_{\tbox{BR}}) \right> \approx  \sqrt{n_1^2 + n_2^2} (n_1n_2)^2 p^4 ,
\end{equation}
\begin{equation}
\label{mSOBGp}
\left< ^mSO(G_{\tbox{BR}}) \right> \approx  \frac{n_1n_2}{\sqrt{n_1^2 + n_2^2}} ,
\end{equation}
\begin{equation}
\label{BSOBGp}
\left< BSO(G_{\tbox{BR}}) \right> \approx  \sqrt{n_1^2 + n_2^2} ,
\end{equation}
\begin{equation}
\label{SOaBGp}
\left< SO_\alpha(G_{\tbox{BR}}) \right> \approx  \left(n_1^\alpha + n_2^\alpha\right)^{1/\alpha} (n_1n_2)^2 p^4 ,
\end{equation}
\begin{equation}
\label{KAabBGp}
\left< KA^1_{\alpha,\beta}(G_{\tbox{BR}}) \right> \approx  n_1n_2p \left[(n_1p)^\alpha + (n_2p)^\alpha\right]^\beta .
\end{equation}
Above we used $|E(G_{\tbox{BR}})|=n_1n_2p$.
In Figs.~\ref{Fig05}(a)-\ref{Fig05}(c), we show that Eqs.~(\ref{SOBGp}-\ref{BSOBGp}) (black-dashed 
lines) indeed describe well the data (thick full curves) for large enough $p$.

As for ER graphs, here for BR graphs the average modified Sombor index and the average first Banhatti-Sombor 
index do not depend on the probability $p$ in the dense limit, see Eqs.~(\ref{mSOBGp},\ref{BSOBGp}).
Also, by recognizing the average degrees $\left< k_{1,2} \right>$ in Eqs.~(\ref{SOBGp},\ref{SOaBGp},\ref{KAabBGp}),
they can be rewritten as
\begin{equation}
\label{SOBGk}
\left< SO(G_{\tbox{BR}}) \right> \approx  \sqrt{n_1^2 + n_2^2} \left( \left< k_1 \right> \left< k_2 \right> \right)^2 ,
\end{equation}
\begin{equation}
\label{SOaBGk}
\left< SO_\alpha(G_{\tbox{BR}}) \right> \approx  \left(n_1^\alpha + n_2^\alpha\right)^{1/\alpha} \left( \left< k_1 \right> \left< k_2 \right>\right)^2 ,
\end{equation}
\begin{equation}
\label{KAabBGk}
\left< KA^1_{\alpha,\beta}(G_{\tbox{BR}}) \right> \approx  |E(G_{\tbox{BR}})| \left(\left< k_1 \right>^\alpha + \left< k_2 \right>^\alpha\right)^\beta .
\end{equation}

Therefore, by plotting $\left< \overline{SO}(G_{\tbox{BR}}) \right>$ vs.~$\left< k_1 \right>\left< k_2 \right>$, 
$\left< ^m\overline{SO}(G_{\tbox{BR}}) \right>$ vs.~$p$, and 
$\left< \overline{BSO}(G_{\tbox{BR}}) \right>$ vs.~$p$ [with 
$\left< \overline{SO}(G_{\tbox{BR}}) \right>=\left< SO(G_{\tbox{BR}}) \right>/\sqrt{n_1^2+n_2^2}$,
$\left< ^m\overline{SO}(G_{\tbox{BR}}) \right>=\sqrt{n_1^2+n_2^2} \left< ^mSO(G_{\tbox{BR}}) \right>/(n_1n_2)$, and
$\left< \overline{BSO}(G_{\tbox{BR}}) \right>=\left< BSO(G_{\tbox{BR}}) \right>/\sqrt{n_1^2+n_2^2}$],
see Figs.~\ref{Fig05}(d-f), we confirm that the curves of these average Sombor indices on BR graphs coincide
in the dense limit, as predicted by Eqs.~(\ref{SOBGk}), (\ref{mSOBGp}) and (\ref{BSOBGp}), respectively.

It is relevant to stress that, while the curves $\left< ^m\overline{SO}(G_{\tbox{BR}}) \right>$ vs.~$p$ and 
$\left< \overline{BSO}(G_{\tbox{BR}}) \right>$ vs~$p$ are properly normalized on the vertical axis, they are 
still not scaled on the $p$-axis. That is, the curves of the main panels in Figs.~\ref{Fig05}(e,f) do not coincide.
However, through a standard scaling analysis (not shown here), 
it is possible to find the scaling parameter $p^*$ 
such that the curves $\left< ^m\overline{SO}(G_{\tbox{BR}}) \right>$ vs.~$\overline{p}$ and 
$\left< \overline{BSO}(G_{\tbox{BR}}) \right>$ vs~$\overline{p}$, with $\overline{p}\equiv p/p^*$, 
fall on top of each other.
Indeed, we found that $p^*=1/\sqrt{n_1^2+n_2^2}$ for $\left< ^m\overline{SO}(G_{\tbox{BR}}) \right>$
and $p^*=\sqrt{n_1^2+n_2^2}/(n_1n_2)$ for $\left< \overline{BSO}(G_{\tbox{BR}}) \right>$.
Thus, as can be seen in the instes of Figs.~\ref{Fig05}(e,f), the curves of the main panels are now
properly scaled when plotted as a function of $\overline{p}$.

It is remarkable to notice that in the case of $n_1=n_2=n/2$, where $\left< k_1 \right>=\left< k_2 \right>=\left< k \right>=np/2$,
we get exactly the same expressions listed in Eqs.~(\ref{SOk}-\ref{KAabk}).
This is verified in Figs.~\ref{Fig06}(d-f) where we plot average Sombor indices
on RG graphs, normalized to $n$, as a function of $\left< k \right>$ showing that all curves are properly scaled.

\section{General scaling of Sombor indices on random graphs}

In the previous Section we have shown that the average value of Sombor indices, normalized to the graph size, 
scale with the average degree $\bra k \ket$ of the corresponding random graph models; we note that this also 
applies to BR graphs when $n_1=n_2$. 
This means that $\bra k \ket$ fixes the average value of any Sombor index for different combinations
of graph parameters; i.e.~the relevant parameter of the random graph models we study here is $\bra k \ket$ 
and not the specific values of the model parameters.
This result highlights the relevance of $\bra k \ket$ in random graph studies.
Moreover, the applicability of Eqs.~(\ref{SOk}-\ref{KAabk}) to the three random graph models we study here
allow us to relate the average value of a given Sombor index $X$ of the three random graph models as
\begin{equation}
\label{XofG}
\frac{\left< X(G_{\tbox{ER}}) \right>}{n} \approx \frac{\left< X(G_{\tbox{RG}}) \right>}{n} 
\approx \frac{\left< X(G_{\tbox{BR}}) \right>}{n}
\quad \quad \mbox{if} \quad \quad 
\left< k_{\tbox{ER}} \right>\approx \left< k_{\tbox{RG}} \right> \approx \left< k_{\tbox{BR}} \right>,
\end{equation}
where $\left< k_{\tbox{ER}} \right>$, $\left< k_{\tbox{RG}} \right>$, and $\left< k_{\tbox{BR}} \right>$
are given in Eqs.~(\ref{k}), (\ref{kRG}), and (\ref{kBR}), respectively.

\begin{figure}[t!]
\begin{center} 
\includegraphics[width=0.99\textwidth]{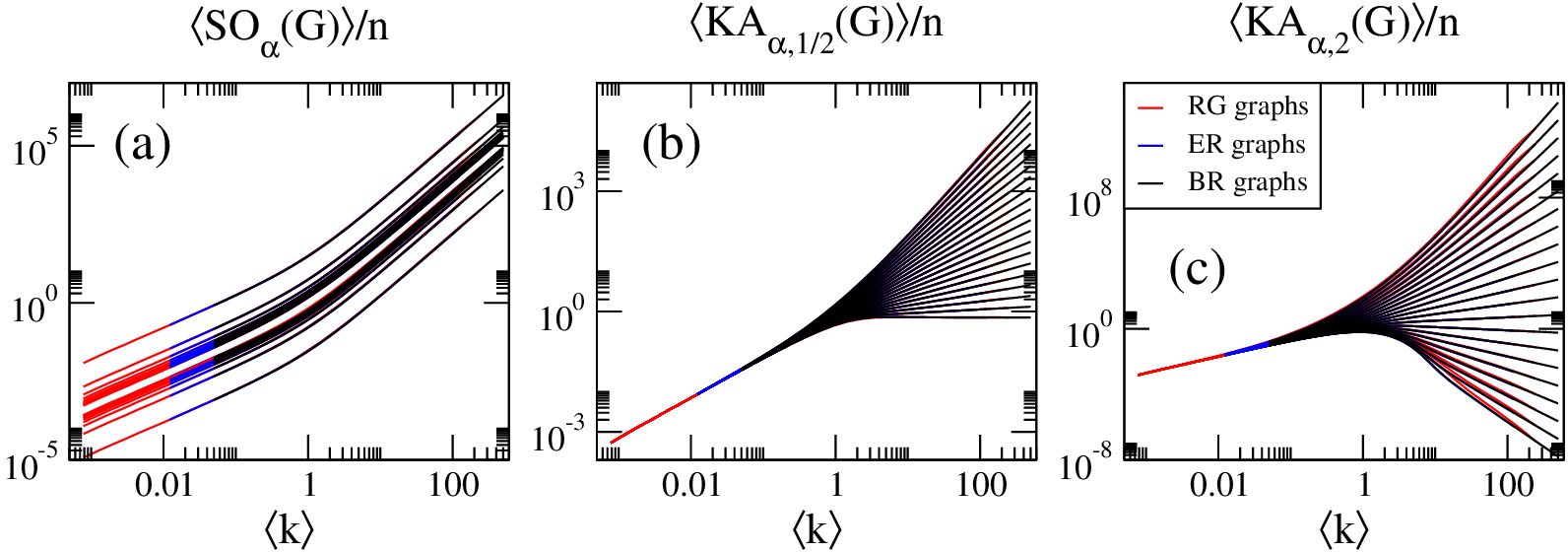}
\caption{\footnotesize{
(a) $\left< SO_\alpha(G) \right>/n$, 
(b) $\left< KA_{\alpha,1/2}(G) \right>/n$, and
(c) $\left< KA_{\alpha,2}(G) \right>/n$
as a function of the average degree $\left< k \right>$ for RG, ER, and BR graphs.
In all panels we show curves for $\alpha\in[-2,2]$ in steps of 0.2 (from bottom to top).
}}
\label{Fig07}
\end{center}
\end{figure}

Now, to verify Eq.~(\ref{XofG}), in Fig.~\ref{Fig07} we compare normalized Sombor indices, 
$\left< X(G) \right>/n$, for ER, RG, and BR graphs, as a function of the corresponding $\left< k \right>$.
Note that to really put Eq.~(\ref{XofG}) to test, we are using graphs of different sizes. 
Indeed, we observe that Eq.~(\ref{XofG}) is satisfied
to a good numerical accuracy; that is, we observe the coincidence of the curves 
$\left< X(G) \right>/n$ vs.~$\left< k \right>$ corresponding to different graphs models.

\section{Sombor indices as complexity mesures for random graphs}

Additionally, we want to recall that in complex systems research there is a continuous search of measures 
that could serve as complexity indicators. In particular, random matrix theory (RMT) has provided us 
with a number of measures able to distinguish between (i) integrable and chaotic (i.e.~non-integrable) 
and (ii) ordered and disordered quantum systems~\cite{M04,H10}.
Such measures are computed from the eigenvalues and eigenvectors of quantum Hamiltonian
matrices. Examples of eigenvalue-based measures are the distribution of consecutive eigenvalue 
spacings, the spectrum rigidity and the ratios between consecutive eigenvalue spacings; while the 
inverse participation ratios and Shannon entropies are popular eigenvector-based complexity 
measures~\cite{M04,H10}.
It is interesting to notice that all these RMT measures have also been successfully applied to study 
networks and graphs since they can be computed from the eigenvalues and eigenvectors of 
adjacency matrices; see e.g.~\cite{MAMRP15,AMGM18,TFM20} and the references therein. Therefore, 
these measures are able to distinguish between graphs composed by mostly isolated vertices and 
mostly connected graphs. Also, through scaling studies of RMT measures it has been possible to 
locate the percolation transition point of random graphs models~\cite{MAMRP15,AMGM18}.
It is worth mentioning that the scaling study of average Sombor indices performed in this paper has 
followed a {\it statistical} RMT approach; that is, from a detailed computational study we have been able to 
identify the average degree as the universal parameter of our random graph models: i.e.~the 
parameter that fixes the average values of the Sombor indices.

Moreover, recently, it has been shown for RG graphs that there is a a huge correlation between the 
average-scaled Shannon entropy (of the adjacency matrix eigenvectors) and two average-scaled 
topological indices~\cite{AMRS20}: the Randi\'c index $R(G)$ and the harmonic index $H(G)$.
We believe that this is a remarkable result because it validates the use of average topological indices
as RMT complexity measures; already suggested in Refs.~\cite{MMRS20,MMRS21} for ER random 
networks.
Now, it is important to stress that not every index could be used as a complexity measure. From our
experience, we conclude
that good candidates should fulfill a particular requirement: they should get well defined values in the 
trivial regimes (just as RMT measures are).
For example, a useful complexity measure for random graphs should be close to zero in the regime of 
mostly isolated vertices while it should become constant above the percolation transition. 
Indeed, this is a property that both $\left< R(G) \right>$ and $\left< H(G) \right>$ have: 
$\left< R(G) \right> \approx \left< H(G) \right> \approx 0$ for mostly isolated vertices
while $\left< R(G) \right>/n \approx \left< H(G) \right>/n \approx 1/2$ once the network is well above 
the percolation transition.

\begin{figure}[t!]
\begin{center} 
\includegraphics[width=0.99\textwidth]{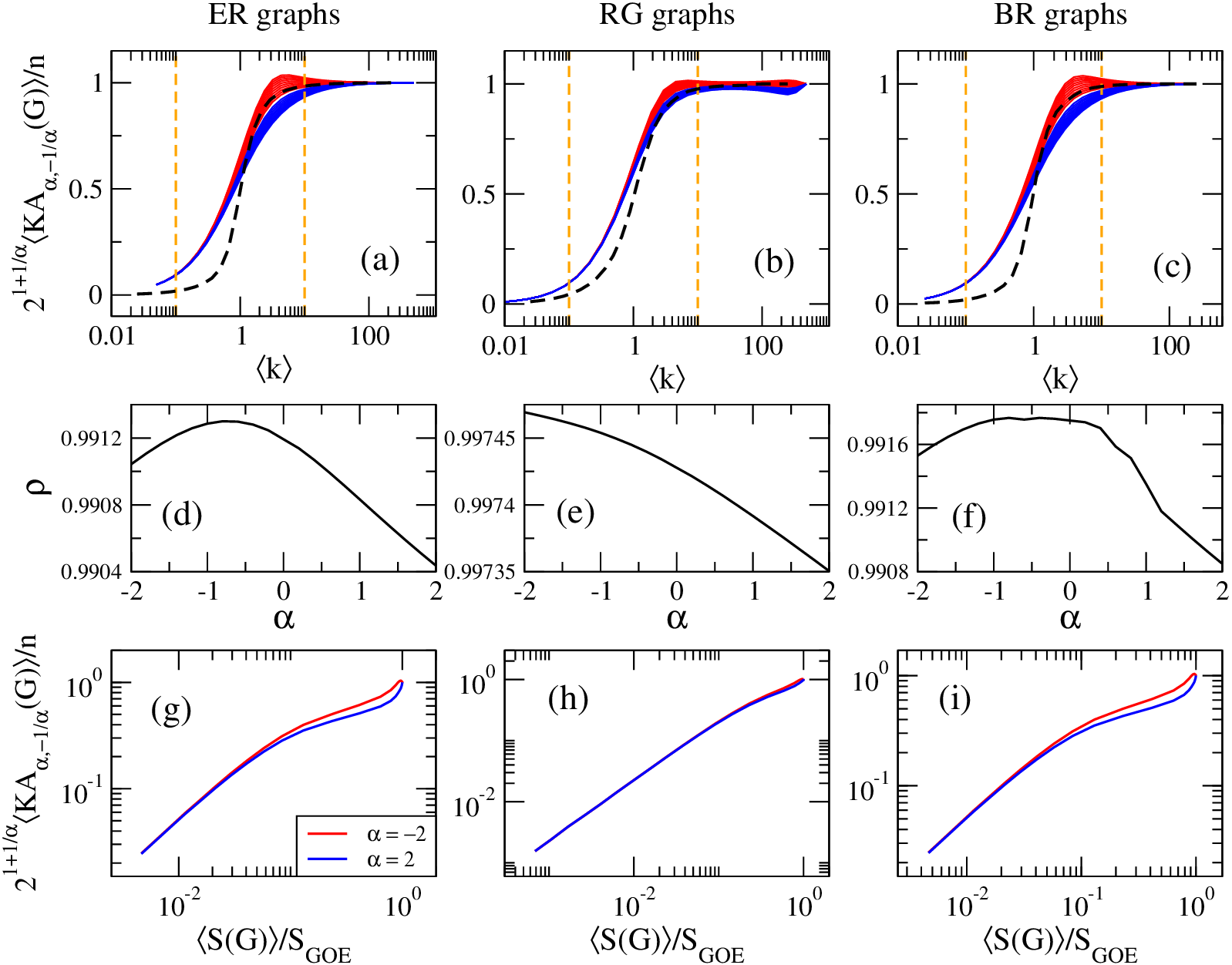}
\caption{\footnotesize{
$2^{1+1/\alpha}\left< KA^1_{\alpha,-1/\alpha}(G) \right>/n$
as a function of the average degree $\left< k \right>$ for 
(a) ER graphs of size $n=500$, (b) RG graphs of size $n=500$, and (c) BR graphs with $n_1=n_2=250$.
In all panels we show curves for $\alpha\in[-2,2]$ in steps of 0.2; except for $\alpha=0$.
Red (blue) lines correspond to $\alpha<0$ ($\alpha>0$).
Vertical orange dashed-lines mark $\left< k \right>=1/10$ and $\left< k \right>=10$, see the text.
Black-dashed lines in (a-c) correspond to the normalized Shannon entropies $\left< S(G) \right>/S_{\tbox{GOE}}$.
(d-f) Pearson's correlation coefficient $\rho$ between $2^{1+1/\alpha}\left< KA^1_{\alpha,-1/\alpha}(G) \right>/n$ 
and $\left< S(G) \right>/S_{\tbox{GOE}}$ as a function of $\alpha$.  
(g-i) Scatter plots of $2^{1+1/\alpha}\left< KA^1_{\alpha,-1/\alpha}(G) \right>/n$ 
vs.~$\left< S(G) \right>/S_{\tbox{GOE}}$ for $\alpha=-2$ and 2.
}}
\label{Fig08}
\end{center}
\end{figure}

Therefore, a straightforward application of our study on Sombor indices is the 
identification of specific Sombor indices as complexity measure candidates. Recall that we particularly 
require, for an average-scaled Sombor index to work as complexity measure, that 
$\left< X(G) \right>/n \approx \mbox{const.}$ for large enough $\left< k \right>$. 
In fact, from Eqs.~(\ref{mSOk}) and~(\ref{BSOk}) we can see that the above condition is fulfilled for 
$\left< ^mSO(G) \right>$ and $\left< BSO(G) \right>$, respectively.
More generally, by properly choosing the values of $\alpha$ and $\beta$ in Eq.~(\ref{KAabk}) we could
also use $\left< KA_{\alpha,\beta}(G) \right>$ as complexity measure. Specifically, for $\beta=-1/\alpha$
we get
\begin{equation}
\label{KAbbk}
\frac{\left< KA^1_{\alpha,-1/\alpha}(G) \right>}{n} \approx \frac{1}{2^{1+1/\alpha}} .
\end{equation}
Note that $\left< KA^1_{\alpha,-1/\alpha}(G) \right>$ reproduces both $\left< ^mSO(G) \right>$ and 
$\left< BSO(G) \right>$ when $\alpha=2$ and $\alpha=-2$, respectively.
Thus, in Fig.~\ref{Fig08} we plot $2^{1+1/\alpha}\left< KA^1_{\alpha,-1/\alpha}(G) \right>/n$ as a function 
of the average degree $\left< k \right>$ for ER, RG, and BR graphs.
From the behavior of the average-scaled indices reported in Fig.~\ref{Fig08} we can identify three
regimes:
(i) a regime of mostly isolated vertices when $\left< k \right> < 1/10$, 
where $2^{1+1/\alpha}\left< KA^1_{\alpha,-1/\alpha}(G) \right>/n \approx 0$, 
(ii) a regime corresponding to mostly connected graphs when $\left< k \right> > 10$, where 
$2^{1+1/\alpha}\left< KA^1_{\alpha,-1/\alpha}(G) \right>/n \approx 1$, and
(iii) a transition regime in the interval $1/10 < \left< k \right> < 10$, which is logarithmically 
symmetric around the percolation transition point $\left< k \right>\approx 1$.
Accordingly, we propose the use of $\left< KA^1_{\alpha,-1/\alpha}(G) \right>$ as 
complexity measure for random graph models.

\subsection{Correlation between the average $KA^1_{\alpha,-1/\alpha}(G)$ index and the average Shannon entropy}

Since we are proposing the use of $\left< KA^1_{\alpha,-1/\alpha}(G) \right>$ as a complexity measure 
for random graphs, it is pertinent to compare it to other standard RMT complexity measure.
To this end we choose the average Shannon entropy $\left< S \right>$ of the adjacency matrix eigenvectors.

In particular we construct randomly weighted adjacency matrices, see e.g.~\cite{AMRS20}, such that we obtain 
well-known RMT ensembles in the limits of: 
(i) isolated vertices (where we get random diagonal adjacency matrices, known in RMT as the Poisson ensemble) and
(ii) complete graphs (where the adjacency matrices become members of the Gaussian Orthogonal Ensemble (GOE)). 
Specifically, for the normalized eigenvector $\Psi^i$, i.e. $\sum_{j=1}^n | \Psi^i_j |^2 =1$, $S$ is
defined as
\begin{equation}
\label{S}
S_i = -\sum_{j=1}^n \left| \Psi^i_j \right|^2 \ln \left| \Psi^i_j \right| ^2 \ .
\end{equation}
Then, we use exact numerical diagonalization to obtain the eigenvectors $\Psi^i$ ($i =1,\ldots,n$) of large 
ensembles of adjacency matrices and compute $\left< S \right>$, where the average is taken over all the 
eigenvectors of all the adjacency matrices of the ensemble.

In Figs.~\ref{Fig08}(a-c) we present $\left< S(G) \right>$, normalized to 
$S_{\tbox{GOE}}\approx \ln (n/2.07)$, for ER, RG and BR graphs; see the black-dashed lines.
From these figures one can observe that $\left< KA^1_{\alpha,-1/\alpha}(G) \right>$ and
$\left< S(G) \right>$ are indeed highly correlated.
To quantify the correlation, in panels Figs.~\ref{Fig08}(d-f) we report the corresponding Pearson's correlation 
coefficient $\rho$, which turns out to be approximately equal to one for all the values of $\alpha$ we consider.
Finally, to validate the high correlation reported by $\rho$, in Figs.~\ref{Fig08}(g-i) we show two examples of 
scatter plots of $2^{1+1/\alpha}\left< KA^1_{\alpha,-1/\alpha}(G) \right>/n$ vs.~$\left< S(G) \right>/S_{\tbox{GOE}}$.

\section{Conclusions}
\label{Conclusions}

In this paper we have performed a thorough computational study of Sombor indices on 
random graphs.
As models of random graphs we have used 
Erd\"os-R\'enyi graphs, random geometric graphs, and bipartite random graphs.

Within a statistical
random matrix theory approach, we show that the average values of Sombor indices, 
normalized to the order of the graph $n$, scale with the graph average degree $\left< k \right>$. 
Thus, we conclude that $\left< k \right>$ is the parameter that fixes the average values of Sombor indices
on random graphs.
Moreover, it is remarkable that we were able to state a scaling law that includes different graph models;
see Eq.~(\ref{XofG}) and Fig.~\ref{Fig07}.

Moreover, we discuss the application of Sombor indices as complexity measures of random graphs and, 
as a consequence, we show that the average first $(\alpha,\beta)-KA$ index (with $\beta=-1/\alpha$), 
normalized to $n$, is 
highly correlated with the averaged-scaled Shannon entropy of the eigenvectors of the graph adjacency 
matrix. That is, $\left< KA^1_{\alpha,-1/\alpha}(G) \right>/n$ may serve as 
complexity measure for random graph models.

We hope that our work may motivate further analytical as well as computational studies of 
Sombor indices on random graphs.

\section*{\bf ACKNOWLEDGEMENTS}
The research of J.M.R. and J.M.S. was supported by a grant from Agencia Estatal de Investigaci\'on (PID2019-106433GBI00/AEI/10.13039/501100011033), Spain.
J.M.R. was supported by the Madrid Government (Comunidad de Madrid-Spain) under the Multiannual Agreement with UC3M in the line of Excellence of University Professors (EPUC3M23), and in the context of the V PRICIT (Regional Programme of Research and Technological Innovation).


\end{document}